\documentstyle[12pt]{article}

\newfont{\myfont}{cmbx10}

\title{ Gauss-type Quadrature Rules for Rational Functions\thanks{
Work supported in part by the National Science Foundation under
grant DMS--9023403.}}

\smallskip
\author{Walter Gautschi} 
\date{\ }

\begin{document}

\maketitle

\begin{quotation}
\noindent
{\small {\bf Abstract}.
When integrating functions that have poles outside the interval of
integration, but are regular otherwise, it is suggested that the
quadrature rule in question ought to integrate exactly
not only polynomials (if any), but also suitable rational
functions. The latter are to be chosen so as to match the most important
poles of the integrand. We describe two methods for generating such
quadrature rules numerically and report on computational experience with
them.}
\end{quotation}
\vspace{.25in}

\noindent
{\bf Introduction}
\bigskip

Traditionally, Gauss quadrature rules are designed to integrate exactly
polynomials of maximum possible degree. This is meaningful for integrand
functions that are ``polynomial-like''. For integrands having poles
(outside the interval of integration) it would be more natural to 
include also rational functions among the functions to be exactly 
integrated. In this paper we consider $n$-point quadrature rules that
exactly integrate $m$ rational functions (with prescribed location and
multiplicity of the poles) as well as polynomials of degree $2n-m-1$,
where $0 \leq m \leq 2n$. The limit case $m = 2n$, in which only
rational functions are being integrated exactly, is a rational
counterpart of the classical Gauss formula; the latter corresponds to
the other limit case $m = 0$.

In \S1 we characterize these new quadrature rules in terms of classical
(polynomial) Gauss formulae with modified weight functions. We also
identify special choices of poles that are of interest in applications.
The computation of the quadrature rules is discussed in \S2, and
numerical examples are given in \S3.

\newpage

\noindent
{\bf 1.  Gauss quadrature for rational functions}
\bigskip

Let $d \lambda$ be a measure on the real line having finite moments
of all orders.
Let $\zeta_\mu \in \mbox{\myfont C}$, $\mu = 1, 2,\ldots, M$,
be distinct real or complex numbers such that
$$
\zeta_\mu \neq 0 ~~ \mbox{and} ~~ 1 + \zeta_\mu t \neq 0 ~~ \mbox{for} ~~
t \in \overline{\mbox{supp} (d \lambda )} , ~~~~ \mu = 1, 2,\ldots, M ~.
\eqno(1.1)
$$
For given integers $m$, $n$ with $1 \leq m \leq 2n$, we wish to find an
$n$-point quadrature rule that integrates exactly (against the measure
$d \lambda$) polynomials of degree $2n-m-1$ as well as the $m$ rational
functions
$$
( 1 + \zeta_\mu t )^{-s} , ~~~~ \mu = 1, 2,\ldots, M , ~~~~
s = 1, 2,\ldots, s_\mu ,
\eqno(1.2)
$$
where $s_\mu \geq 1$ and
$$
{\displaystyle \sum_{\mu=1}^M s_\mu = m } .
\eqno(1.3)
$$
In the extreme case $m=2n$ (where polynomials of degree --1 are
understood to be identically zero) the formula integrates exactly $2n$
rational functions (with poles of multiplicities $s_\mu$ at 
$-1/\zeta_\mu$), but no nontrivial polynomials. The formula, therefore,
can be thought of as the rational analogue of the classical Gauss
formula; the latter corresponds to the other limit case $m = M = 0$.

The solution of our problem is given by the following theorem.
\bigskip

THEOREM 1.1.
{\it Define
$$
\omega_m (t) = \prod_{\mu =1}^M (1 + \zeta_\mu t)^{s_{\mu}} , 
\eqno(1.4)
$$
a polynomial of degree m.
Assume that the measure $d \lambda / \omega_m$ admits a} (polynomial)
{\it n-point Gaussian quadrature formula
$$
\int_{\mbox{\myfont R}} f(t) ~
\frac{d \lambda (t)}{\omega_m (t)} ~=~
\sum_{\nu =1}^n w_\nu ^G f (t_\nu ^G ) 
+ R_n ^G (f), ~~~~ R_n ^G ( \mbox{\myfont P}_{2n-1}) = 0,
\eqno(1.5)
$$
with nodes $t_\nu^G$ contained in the support of $d \lambda$,
$$
t_\nu^G  \in  \mbox{\rm supp} ( d \lambda ).
\eqno(1.6)
$$
Define
$$
t_\nu = t_\nu ^G , ~~~~ \lambda_\nu = w_\nu ^G  \omega_m (t_\nu ^G ), ~~~~
\nu = 1, 2,\ldots, n .
\eqno(1.7)
$$
Then
$$
\int_{\mbox{\myfont R}} g(t) d \lambda (t) = \sum_{\nu=1}^n \lambda_\nu
g (t_\nu ) + R_n (g) ,
\eqno(1.8)
$$
where}
$$
\renewcommand{\arraystretch}{2}
R_n (g) = 0 ~~ i f ~~ \left\{ \begin{array}{l}
g(t)  = (1 + \zeta_\mu t)^{-s} , ~~~
\mu =1,2,\ldots, M; ~
s = 1, 2,\ldots, s_\mu, \\
g \in \mbox{\myfont P}_{2n-m-1} ~.
\end{array} \right.
\eqno(1.9)
$$
\renewcommand{\arraystretch}{1}

{\it Conversely}, (1.8) {\it with} $t_\nu \in 
\mbox{\rm supp} ( d \lambda )$ {\it and} (1.9) {\it imply} (1.5), (1.6) {\it with $t_\nu^G ,
w_\nu^G$ as defined in} (1.7).
\bigskip

{\it Remark}. Theorem 1.1, for real $\zeta_\mu$ and either all $s_\mu = 1$ and
$m = 2n$, or all but one $s_\mu = 2$ and $m = 2n-1$,
is due to Van Assche and Vanherwegen [13]. The quadrature
rule (1.5), especially its convergence properties for analytic
functions $f$, has previously been studied by L\'{o}pez and Ill\'{a}n
[9, 10].
\bigskip

{\it Proof of Theorem} 1.1.
Assume first (1.5), (1.6). For $\mu = 1, 2,\ldots, M; s = 1, 2,\ldots, s_\mu$, define
$$
q_{\mu ,s} (t) ~=~ \frac{\omega_m (t)}{(1 + \zeta_\mu t)^s }~.
\eqno(1.10)
$$
Since $m \leq 2n$ and $s \geq 1$, we have
$q_{\mu ,s} \in \mbox{\myfont P}_{m-s} \subset \mbox{\myfont P}_{2n-1}$,
and therefore, by (1.5),
\renewcommand{\arraystretch}{2}
$$
\begin{array}{c}
{\displaystyle \int_{\mbox{\myfont R}} ~
{\displaystyle \frac{d \lambda (t)}{(1 + \zeta_\mu t)^s} } } ~=~
{\displaystyle \int_{\mbox{\myfont R}} q_{\mu ,s} (t) ~
{\displaystyle 
\frac{d \lambda (t)}{\omega_m (t)} } }~=~
{\displaystyle \sum_{\nu =1}^n w_\nu^G q_{\mu ,s} (t_\nu^G ) } \\
= {\displaystyle \sum_{\nu =1}^n w_\nu^G } ~
{\displaystyle
\frac{\omega_m (t_\nu^G )}{(1 + \zeta_\mu t_\nu^G )^s} } ~=~
{\displaystyle \sum_{\nu =1}^n }~
{\displaystyle
\frac{\lambda_\nu}{(1 + \zeta_\mu t_\nu )^s} } ~,
\end{array}
$$
\renewcommand{\arraystretch}{1}where (1.7)
has been used in the last step and none of the denominators on the
far right vanishes by (1.6) and (1.1).
This proves the assertion in the top line of (1.9).
The bottom part of (1.9) follows similarly: Let $p$ be an arbitrary
polynomial in $\mbox{\myfont P}_{2n-m-1}$.
Then, since $p \, \omega_m \in \mbox{\myfont P}_{2n-1}$, again by
(1.5) and (1.7),
\renewcommand{\arraystretch}{2}
$$
\begin{array}{c}
{\displaystyle
\int_{\mbox{\myfont R}} p(t) d \lambda (t) = \int_{\mbox{\myfont
R}} p(t) \omega_m (t) \frac{d \lambda (t)}{\omega_m (t)} } \\
{\displaystyle
= \sum_{\nu =1}^n w_\nu^G p (t_\nu^G ) \omega_m (t_\nu^G ) =
\sum_{\nu =1}^n \lambda_\nu p(t_\nu )} . 
\end{array}
$$
\renewcommand{\arraystretch}{1}

To prove the converse, we first note that $w_\nu^G$ is well
defined by (1.7), since $\omega _m (t_\nu ) \neq 0$ by the assumption
on $t_\nu$ and (1.1). One then easily verifies that (1.5) holds for
all polynomials (1.10) (of degree $~ < m$) and all polynomials of the
form $p \, \omega_m$ where $p \in \mbox{\myfont P}_{2n-1-m}$. The
collection of these polynomials, however, spans 
$\mbox{\myfont P}_{2n-1}. ~~~~ \Box$
\bigskip

We will concentrate on six special choices of the parameters
$\zeta_\mu$ that are of interest in applications.
\bigskip

{\it Case} 1 (Simple real poles).
All $s_\mu = 1$ in (1.2) (hence $M = m$), and all $\zeta_\mu$ are real,
distinct, and nonzero,
$$
\zeta_\nu = \xi_\nu \in \mbox{\myfont R} , ~~~~
\xi_\nu \neq 0 , ~~~~
\nu = 1, 2,\ldots, m .
\eqno(1.11\mbox{a})
$$
In this case the polynomial $\omega_m$ has the form
$$
\omega_m (t) = \prod_{\nu =1}^m (1 + \xi_\nu t) , ~~~~
\xi_\nu \in \mbox{\myfont R} .
\eqno(1.11\mbox{b})
$$
If the support of $d \lambda$ is an interval, $\omega_m$ does not
change sign on it because of (1.1).
\bigskip

{\it Case} 2 (Simple conjugate complex poles).
All $s_\mu = 1$ (hence $M = m$), $m$ even, and the $\zeta_\mu$ occur in 
$m/2$ (distinct) pairs of
conjugate complex numbers (cf. [9]),
$$
\zeta_\nu = \xi_\nu + i \eta_\nu , ~~~~
\zeta_{\nu + m/2} = \xi_\nu - i \eta_\nu , ~~~~
\nu =  1, 2,\ldots, m/2 ,
\eqno(1.12\mbox{a})
$$
where $\xi_\nu \in \mbox{\myfont R}$ and
$\eta_\nu > 0$.
Here,
$$
\omega_m (t) = \prod_{\nu =1}^{m/2} [ ( 1 + \xi_\nu t)^2 +
\eta_\nu^2 t^2 ] ,
\eqno(1.12\mbox{b})
$$
which is strictly positive for real $t$.
\bigskip

{\it Case} 2$^\prime$ (Simple conjugate complex poles plus a simple real pole).
All $s_\mu =1$ (hence $M = m$), $m$ (odd) $\geq 3$, and,
slightly changing the indexing of the $\zeta$'s,
$$
\zeta_0 \in \mbox{\myfont R} , ~~~
\zeta_\nu = \xi_\nu + i \eta_\nu , ~~~
\zeta_{\nu + (m-1)/2} = \xi_\nu - i \eta_\nu , ~~~
\nu = 1, 2 ,\ldots, (m-1)/2 , 
\eqno(1.13\mbox{a})
$$
where $\zeta_0 = \xi_0 \neq 0$ and
$\xi_\nu \in \mbox{\myfont R}, ~ \eta_\nu > 0$ for $1 \leq \nu 
\leq (m-1)/2$.
Then
$$
\omega_m (t) = (1 + \xi_0 t) 
\prod_{\nu =1}^{(m-1)/2} [(1 + \xi_\nu t)^2 + \eta_\nu^2 t^2 ] .
\eqno(1.13\mbox{b})
$$
\bigskip

{\it Case} 3 (Real poles of order 2).
All $s_\mu = 2$ in (1.2) (hence $2M = m$), and all $\zeta_\mu$ are
nonzero, real and distinct,
$$
\zeta_\nu = \xi_\nu \in \mbox{\myfont R} , ~~~~
\xi_\nu \neq 0 , ~~~~
s_\nu = 2, ~~~~
\nu = 1, 2,\ldots, m/2 .
\eqno(1.14\mbox{a})
$$
The polynomial $\omega_m$ now has the form
$$
\omega_m (t) = \prod_{\nu =1}^{m/2} (1 + \xi_\nu t)^2
\eqno(1.14\mbox{b})
$$
and is nonnegative for real $t$, and positive on the support of $d
\lambda$.
\bigskip

{\it Case} 3$^\prime$ (Real poles of order 2 plus a simple real pole).
Here, all $\zeta_\mu = \xi_\mu$
are nonzero, real and distinct,
$s_\nu = 2$ for $\nu = 1, 2,\ldots, M-1$ and $s_M = 1$.
Thus, $m = 2M-1$, and 
$$
{\displaystyle
\omega_m (t) = (1 + \xi_M t) \prod_{\nu =1}^{M-1} (1 + \xi_\nu t)^2 , ~~
M = \frac{m+1}{2} , ~~
m \mbox{(odd)} \geq 3 . }
\eqno(1.15)
$$
If $M = n$, i.e., $m = 2n-1$, the quadrature rule (1.5) is then
identical with the ``orthogonal quadrature rule'' of [13], having
as nodes the zeros of the rational function
$(1 + \zeta_n t)^{-1} + \sum_{\nu =1}^{n-1} a_\nu (1 + \zeta_\nu t
)^{-1}$ which is orthogonal (relative to the measure $d \lambda$)
to 1 and to $(1 + \zeta_\mu t)^{-1}$,
$\mu = 1, 2,\ldots, n-1$.
As in Case 1, the polynomial $\omega_m$ preserves its sign on the interval on
which $d \lambda$ is supported.
\bigskip

{\it Case} 4 (Conjugate complex poles of order 2).
All $s_\mu = 2$ (hence $2M = m$), $m = 0$ (mod 4), and the $\zeta_\mu$,
$\mu = 1, 2,\ldots, m/2$, occur in $m/4$ conjugate complex pairs,
similarly as in Case 2.
Thus,
$$
\omega_m (t) = \prod_{\nu =1}^{m/4} [( 1+ \xi_\nu t)^2 +
\eta_\nu^2 t^2]^2 .
\eqno(1.16)
$$

In all six cases, the measure $d \lambda / \omega_m$ admits a
Gaussian $n$-point formula for each $n = 1, 2, 3,\ldots,$ so that
the assumption of Theorem 1.1 is fulfilled for each $n$.

Putting $f(t) = \omega_m (t) g(t)$ in (1.5) and using (1.7), we
get
$$
\int_{\mbox{\myfont R}} g(t) d \lambda (t) = \sum_{\nu =1}^n \lambda_\nu
g (t_\nu ) + R_n ^G ( \omega_m g ),
\eqno(1.17)
$$
where from the well-known expression for the remainder term of
Gaussian quadrature rules, one has
$$
R_n ^G ( \omega_m g)  = \gamma_n ( \omega_m g )^{(2n)} ( \tau
), ~~~~
\gamma_n ~=~ \frac{\hat{\beta}_0 \hat{\beta}_1 \cdots
\hat{\beta}_n}{(2n) !} ~.
\eqno(1.18)
$$
Here, $\tau$ is some number in the smallest interval containing the
support of $d \lambda$, and
$\hat{\beta}_k = \beta_k (d \lambda / \omega_m )$ are the
$\beta$-recursion coefficients for the measure $d \lambda /
\omega_m$ (cf. (2.1) below).
The latter are computed as part of the algorithms to be described in
the next section.

\bigskip

\medskip

\noindent
{\bf 2.  Computation of the quadrature rule (1.5)}
\bigskip

We propose essentially two methods for generating the basic
quadrature rule (1.5), the first being most appropriate if the support of $d
\lambda$ is a finite interval, the other more effective, though
possibly slower, when the support interval of $d \lambda$ is
unbounded.

\bigskip

{\bf 2.1.  Method based on partial fraction decomposition and
modification algorithms}.
To compute the $n$-point
formula (1.5), it suffices to compute the $n$th-degree orthogonal
polynomial $\hat{\pi}_n ( ~\cdot~ ) = \pi_n (~\cdot ~; d \hat{\lambda
})$ relative to the measure $d \hat{\lambda} = d \lambda /
\omega_m$, or, more precisely, the recursion coefficients
$\hat{\alpha}_k = \alpha_k ( d \hat{\lambda} )$,
$\hat{\beta}_k = \beta_k (d \hat{\lambda})$,
$k = 0, 1,\ldots,$ $n-1$, in the three-term recurrence relation
satisfied by these (monic) polynomials:
\renewcommand{\arraystretch}{2}
$$
\begin{array}{l}
\hat{\pi}_{k+1} (t) = (t - \hat{\alpha}_k )
\hat{\pi}_k (t) - \hat{\beta}_k \hat{\pi}_{k-1} (t) , \\
\hspace{1.4in}k = 0, 1,\ldots, n-1 , \\
\hat{\pi}_0 (t) = 1 , ~~~~ \hat{\pi}_{-1} (t) = 0 .
\end{array}
\eqno(2.1)
$$
\renewcommand{\arraystretch}{1}The nodes
$t_\nu ^G$ and weights $w_\nu ^G$ in (1.5) can then be obtained by standard
techniques via an eigensystem problem for the (symmetric, tridiagonal)
Jacobi matrix of order $n$
having the $\hat{\alpha}_k$,
$k = 0, 1,\ldots, n-1$, on the diagonal, and
$\sqrt{\hat{\beta}_k}$, $k = 1, 2,\ldots, n-1$,
on the side diagonals (see, e.g., [7], [4, \S6]).
The coefficients $\hat{\alpha}_k$, $\hat{\beta}_k$ in turn are
expressible in terms of the orthogonal polynomials $\hat{\pi}_k$ as
\renewcommand{\arraystretch}{2}
$$
\begin{array}{c}
\hat{\alpha}_k ~=~{\displaystyle \frac{(t \hat{\pi}_k , \hat{\pi}_k )}
{( \hat{\pi}_k , \hat{\pi}_k )}} ~, ~~ 0 \leq k \leq n-1 , \\
\hat{\beta}_0 = ( \hat{\pi}_0 , \hat{\pi}_0 ) , ~
\hat{\beta}_k ~=~
{\displaystyle \frac{(\hat{\pi}_k , \hat{\pi}_k )}{( \hat{\pi}_{k-1} ,
\hat{\pi}_{k-1} )}} ~,  ~~ 1 \leq k \leq n-1 ,
\end{array}
\eqno(2.2)
$$
\renewcommand{\arraystretch}{1}where
$(~\cdot ~,~ \cdot ~)$ denotes the inner product
$$
(u, \nu ) = \int_{\mbox{\myfont R}} u(t) \nu (t) d \hat{\lambda}
(t) .
\eqno(2.3)
$$
(If the error constant
$\gamma_n$ in (1.18) is desired, one needs to compute, in addition,
$\hat{\beta}_n$.)

The basic idea of computing the coefficients in (2.2) is as
follows.
Suppose we can construct an $N$-point
quadrature rule for $d \hat{\lambda} = d \lambda / \omega_m$,
where $N > n$, which is exact for polynomials of degree $\leq 2n-1$:
$$
\int_{\mbox{\myfont R}} p(t) d \hat{\lambda} (t) = \sum_{k=1}^N W_k
p (T_k ) , ~~~~
p \in \mbox{\myfont P}_{2n-1} .
\eqno(2.4)
$$
Here the weights $W_k$ are not necessarily all positive.
Denote the discrete measure implied by the sum on the right by
$d \Lambda_N$:
$$
\int_{\mbox{\myfont R}} p(t) d \Lambda_N (t) = \sum_{k=1}^N W_k
p(T_k ) .
\eqno(2.5)
$$
Then from the formulae in (2.2) one easily sees by induction that
\renewcommand{\arraystretch}{2}
$$
\begin{array}{l}
\alpha_k (d \hat{\lambda} ) = \alpha_k ( d \Lambda_N ) , ~~~~ \\
\beta_k (d \hat{\lambda} ) = \beta_k ( d \Lambda_N ) , ~~~~
\end{array}
k = 0, 1,\ldots, n-1 .
\eqno(2.6)
$$
\renewcommand{\arraystretch}{1}Thus,
the desired recursion coefficients are the first $n$ of the $\alpha$-
and $\beta$-coefficients belonging to the discrete measure $d
\Lambda_N$.
These can be generated by Stieltjes's procedure, which is
implemented in the routine {\tt sti} of [4].
(The faster routine {\tt lancz} of [4], implementing the Lanczos
method, would also be applicable here, even though $d \Lambda_N$ is not
necessarily a positive measure.)

We next show how a quadrature rule of type (2.4), with
$N = O(mn)$, can be constructed by means of partial fraction
decomposition and suitable modification algorithms.
For this, we consider separately Cases 1--3$^\prime$ identified
in \S1.
The analysis of Case 4 becomes so tedious that we will not
pursue it any further in this context; see, however, \S2.2.
\bigskip

2.1.1.
{\it Simple real poles}.
We set up the partial fraction decomposition of $1/ \omega_m$ in
the form
$$
\frac{1}{\omega_m (t)} ~=~
\frac{1}{\prod_{\nu =1}^m (1 + \xi_\nu t )} ~=~
\sum_{\nu =1}^m ~
\frac{c_\nu}{t + (1/ \xi_\nu )} ~,
\eqno(2.7)
$$
where
$$
c_\nu ~=~
\frac{\xi_\nu^{m-2}}
{\prod_{\stackrel{\mu=1}{\mu \neq \nu}}^m ( \xi_\nu - \xi_\mu )}, ~~
\nu = 1, 2,\ldots, m ,
\eqno(2.8)
$$
and an empty product in (2.8) (when $m = 1$) is to be taken as 1.
Then, with $d \hat{\lambda} = d \lambda / \omega_m$,
$$
\int_{\mbox{\myfont R}} p(t) d \hat{\lambda} (t) =
\sum_{\nu =1}^m \int_{\mbox{\myfont R}} p(t) ~
\frac{c_\nu d \lambda (t)}{t + (1/ \xi_\nu )} ~.
$$
The integrals on the right involve measures 
$c_\nu d \lambda$ modified by linear divisors.
For such measures, the associated recursion coefficients can be
obtained from those of $c_\nu d \lambda$ (assumed known) by a
suitable modification algorithm (cf. [4, \S5]).
Unless $x_\nu = -1/ \xi_\nu$ is very close to the support interval
of
$d \lambda$, the most appropriate algorithm is the one embodied in
the routine {\tt gchri} of [4] with {\tt iopt} = 1.
Otherwise,
the routine {\tt chri} of [4] (again with {\tt iopt} =
1) is preferable.
A basic ingredient of the routine {\tt gchri} is the modified
Chebyshev algorithm (cf. [2, \S2.4]) using modified moments
$\int_{\mbox{\myfont R}} \pi_k (t;d \lambda ) c_\nu d \lambda (t)
/ (t - x_\nu )$,
$k = 0, 1, 2,\ldots, 2n-1$.
These in turn are generated by backward recurrence as minimal
solution of the three-term recurrence relation for the measure $d
\lambda$; cf. [1, \S5].

Having obtained, in whichever way, the first $n$ of the $\alpha$-
and $\beta$-coefficients for the modified measure
$c_\nu d \lambda (t) / (t - x_\nu )$, and hence the Gaussian
quadrature formula\footnote{In order to produce positive
$\beta$-coefficients, as required in the routine for Gauss
quadrature formulae, one inputs the measure
$| c_\nu / (t-x_\nu ) | d \lambda (t)$ and, if this entails a
change of sign, reverses the sign of all Gauss weights after exiting
from the Gauss quadrature routine.}
$$
\int_{\mbox{\myfont R}} p(t) ~
\frac{c_\nu d \lambda (t)}{t + (1/ \xi_\nu )} =
\sum_{r=1}^n w_r^{( \nu )} p(t_r^{( \nu )} ) , ~~~~
p \in \mbox{\myfont P}_{2n-1} ,
\eqno(2.9)
$$
via eigensystem techniques, we then get
$$
\int_{\mbox{\myfont R}} p(t) ~
\frac{d \lambda (t)}{\omega_m (t)} ~=~
\sum_{\nu =1}^m \int_{\mbox{\myfont R}} p(t) ~
\frac{c_\nu d \lambda (t)}{t + (1/ \xi_\nu )} 
$$
$$
=~ \sum_{\nu =1}^m \sum_{r=1}^n w_r^{( \nu )} p (t_r^{ (\nu )}
) ,
~~~ p \in \mbox{\myfont P}_{2n-1} ,
$$
hence the desired quadrature rule (2.4), with $N = mn$ and
\renewcommand{\arraystretch}{2}
$$
\begin{array}{l}
T_{( \nu -1 )n+r} = t_r^{( \nu )} , ~~~~ \\
W_{( \nu -1 )n+r} = w_r^{( \nu )} , ~~~~ \end{array}
\nu = 1, 2,\ldots, m; ~~ r = 1, 2,\ldots, n . 
\eqno(2.10)
$$
\renewcommand{\arraystretch}{1}

The procedure described works best if the support of $d \lambda$ is
a finite interval. Other\-wise,
the modified Chebyshev algorithm underlying the
procedure is likely to suffer from ill-conditioning; cf. Example 3.4.
Another difficulty that may adversely affect the accuracy of the
results, in particular if $m = 2n$, is the possibility that the constants
$c_\nu ~\mbox{sgn}_{t \in \mbox{\small supp} (d \lambda )} (t+1/ \xi_\nu
)$ become very large and alternate in sign; cf. Example 3.2.
This will cause serious cancellation errors in evaluating inner
products relative to the measure $d \Lambda_N$ (there being blocks of
weights $W_k$ which are very large positive alternating with blocks
of weights which are very large negative).
In such cases, either $m$ has to be lowered, perhaps down to $m = 1$,
or else the method discussed in \S2.2 invoked, which will be more 
effective (but possibly more expensive).
\bigskip

2.1.2.
{\it Simple conjugate complex poles}.
We now consider Case 2 of \S1, i.e., conjugate complex parameters
$\zeta_\nu = \xi_\nu + i \eta_\nu$,
$\zeta_{\nu + m/2} = \bar{\zeta}_\nu$, where $\xi_\nu \in
\mbox{\myfont R}$,
$\eta_\nu > 0$ and $m$ is even.
In this case, an elementary computation yields the partial fraction
decomposition
$$
\frac{1}{\omega_m (t)} ~=~
\sum_{\nu =1}^{m/2} ~
\frac{c_\nu + d_\nu t}{\left( t + \frac{\xi_\nu}{\xi_\nu^2 +
\eta_\nu^2} \right)^2 + \left( \frac{\eta_\nu}{\xi_\nu^2 +
\eta_\nu^2} \right)^2 } ~, ~~~~
t \in \mbox{\myfont R} ,
\eqno(2.11)
$$
where
\renewcommand{\arraystretch}{2}
$$
\begin{array}{l}
c_\nu ~=~ \frac{1}{\eta_\nu} \left(
\frac{\xi_\nu}{\xi_\nu^2 + \eta_\nu^2} ~
\mbox{Im}~ p_\nu ~+~ \frac{\eta_\nu}{\xi_\nu^2 + \eta_\nu^2 } ~
\mbox{Re}~ p_\nu \right) ,\\
d_\nu ~=~ \frac{1}{\eta_\nu} \mbox{Im}~ p_\nu
\end{array} \eqno(2.12)
$$
\renewcommand{\arraystretch}{1}and
$$
p_\nu = \prod_{\stackrel{\mu=1}{\mu \neq \nu}}^{m/2} ~
\frac{( \xi_\nu + i \eta_\nu)^2}
{( \xi_\nu - \xi_\mu )^2 -
(\eta_\nu^2 - \eta_\mu^2 ) +
2 i  \eta_\nu ( \xi_\nu - \xi_\mu )} ~
\eqno(2.13)
$$
with $p_1 = 1$ if $m = 2$. One can then proceed as in \S2.1.1, except
that the modification of
the measure $d \lambda$ now involves multiplication by a
nonconstant linear function (if $d_\nu \neq 0$) in addition to
division by a quadratic.
The former modification is handled by the routine {\tt chri} of
[4] with {\tt iopt} = 1, the latter by the routine {\tt gchri} with
{\tt iopt} = 2 (or, if more appropriate, by {\tt chri} with {\tt
iopt} = 5). The quadrature rule (2.4) so obtained has $N = mn/2$.

If the poles $-1/\zeta_\mu$ are located in conjugate pairs on a line
parallel to the imaginary axis, then by an elementary calculation
one can show that all $p_\nu$ are real, hence $d_\nu = 0$, and there
is no need to call {\tt chri}.
\bigskip

2.1.$2^\prime$.
{\it Simple conjugate complex poles plus a simple real pole}.
We are now in Case 2$^\prime$ of \S1, with $m$ odd,
$\zeta_0 = \xi_0 \in \mbox{\myfont R}$ and the remaining
$\zeta_\mu$ conjugate complex as in Case 2.
This yields
$$
\frac{1}{\omega_m (t)} ~ = ~
\frac{c_0^\prime}{t + (1/ \xi_0 )} ~+
\sum_{\nu =1}^{(m-1)/2} ~
\frac{c_\nu^\prime + d_\nu^\prime t}
{\left( t + \frac{\xi_\nu}{\xi_\nu^2 + \eta_\nu^2} \right)^2 +
\left( \frac{\eta_\nu}{\xi_\nu^2 + \eta_\nu^2} \right)^2 } ~, ~~
t \in \mbox{\myfont R} ,
\eqno(2.14)
$$
where
$$
c_0^\prime = ~\frac{\xi_0^{m-2}}{\prod_{\nu=1}^{(m-1)/2} [( \xi_0 - \xi_\nu )^2 + \eta_\nu^2 ]} ~, 
$$
$$
c_\nu^\prime = ~\frac{1}{\eta_\nu} ~ \left( \frac{\xi_\nu}{\xi_\nu^2 + \eta_\nu^2} ~ \mbox{Im} ~ p_\nu^\prime + ~ \frac{\eta_\nu}{\xi_\nu^2 + \eta_\nu^2} ~\mbox{Re} ~ p_\nu^\prime \right) ,
\eqno(2.15)
$$
$$
d_\nu^\prime  = ~ \frac{1}{\eta_\nu} ~\mbox{Im} ~ p_\nu^\prime
$$
and
$$
p_\nu^\prime = ~
\frac{\xi_\nu + i \eta_\nu}{\xi_\nu - \xi_0 + i \eta_\nu} ~
p_\nu ,
\eqno(2.16)
$$
with $p_\nu$ the same as in (2.13) with $m$ replaced by $m - 1$.
The technique called for is a combination of the one in \S2.1.1,
to deal with the first term in (2.14), and the one in
\S2.1.2, to deal with the remaining terms, and yields a quadrature
rule (2.4) with $N = (m+1)n/2$.

\bigskip
2.1.3.  {\it Real poles of order 2}.
This is Case 3 of \S1, and leads to the partial fraction
decomposition
$$
\frac{1}{\omega_m (t)} ~=~
\sum_{\nu =1}^{m/2} ~
\left( \frac{c_\nu}{t+1/ \xi_\nu} ~+~
\frac{d_\nu}{(t+1/ \xi_\nu )^2 } \right) ,
\eqno(2.17)
$$
$$
c_\nu ~=~ - ~
\frac{2 \xi_\nu^{m-3} \sum_{\stackrel{\mu=1}{\mu \neq \nu}}^{m/2} ~
\frac{\xi_\mu}{\xi_\nu - \xi_\mu} }
{\prod_{\stackrel{\mu=1}{\mu \neq \nu}}^{m/2} ( \xi_\nu - \xi_ \mu )^2}
\eqno(2.18)
$$
$$
d_\nu ~=~
\frac{\xi_\nu^{m-4}}{\prod_{\stackrel{\mu=1}{\mu \neq \nu}}^{m/2}
(\xi_\nu - \xi_\mu )^2} ~,
\eqno(2.19)
$$
where $c_1 = 0$, $d_1  = \xi_1^{-2}$ when $m=2$. Here, $N = mn$ in
(2.4).
\bigskip

2.1.$3^\prime$.  {\it Real poles of order} 2 {\it plus a simple real pole}.
Similarly as in \S2.1.3, the partial fraction decomposition has now
the form
\renewcommand{\arraystretch}{2}
$$
\displaystyle{
\frac{1}{\omega_m (t)} = 
\frac{c_M^\prime}{t + 1/ \xi_M} +
\sum_{\nu =1}^{M-1} \left(
\frac{c_\nu^\prime}{t + 1/ \xi_\nu} +
\frac{d_\nu^\prime}{(t + 1/\xi_\nu )^2 } \right) , ~~
M = (m+1) / 2, ~~ m ~\mbox{odd} , }
\eqno(2.20)
$$
$$
\begin{array}{lll}
c_M^\prime & = & \displaystyle{
\frac{\xi_M^{m-2}}{\prod_{\nu=1}^{M-1} ( \xi_M - \xi_\nu )^2} ~, } \\
c_\nu^\prime & = - &  \displaystyle{
\frac{\xi_\nu^{m-3} \left( \xi_M + 2 ( \xi_\nu - \xi_M )
\sum_{\stackrel{\mu =1}{\mu \neq \nu}}^{M-1} \frac{\xi_\mu}{\xi_\nu
- \xi_\mu} \right) }{(\xi_\nu - \xi_M )^2 \prod_{\stackrel{\mu
=1}{\mu \neq \nu}}^{M-1} (\xi_\nu - \xi_\mu )^2 } ~, } \\
d_\nu^\prime & = & \displaystyle{
\frac{\xi_\nu^{m-4}}{( \xi_\nu - \xi_M ) \prod_{\stackrel{\mu =1}{\mu
\neq \nu }}^{M-1} (\xi_\nu - \xi_\mu )^2 } ~. }
\end{array}
$$
\renewcommand{\arraystretch}{1}Empty
sums and products (when $M = 2$) have their
conventional values 0 and 1, respectively. Again, $N = mn$ in (2.4).

The presence of two terms in the summations of (2.17) and (2.20)
complicates matters considerably, as they call for two
applications of the routine {\tt gchri}:
First, we must generate sufficiently many of the recursion
coefficients for the measure $d \lambda (t) / (t - x_\nu )$,
$x_\nu = - 1/ \xi_\nu$, in order next to generate the desired
recursion coefficients for $d \lambda (t) / (t - x_\nu )^2$ by
backward recursion -- a recursion based on the recurrence relation
generated in the first application of {\tt gchri} (which in turn
requires backward recursion!).
The procedure nevertheless works well if the $x_\nu$ are not too
close to the support interval of $d \lambda$; see Example 3.3.
\bigskip

{\bf 2.2.  Discretization method}.
In this method, the inner product (2.3) is approximated by a
discrete (positive) inner product,
$$
(u,v) = \int_{\mbox{\myfont R}} u(t) v(t) ~\frac{d \lambda (t)}
{\omega_m (t)} ~\approx \sum_{k=1}^N
\omega_k^{(N)} u (\tau_k^{(N)} )
v (\tau_k^{(N)} ) =: (u,v)_N , ~~ N > n,
\eqno(2.21)
$$
whereupon the formulae (2.2) are applied with the inner product
$( ~ \cdot ~, \cdot ~)$ replaced by
$( ~ \cdot ~, \cdot ~)_N$ throughout.
This yields approximations
$$
\hat{\alpha}_{k,N} \approx \hat{\alpha}_k , ~~~~
\hat{\beta}_{k,N} \approx \hat{\beta}_k , ~~~~
k = 0, 1,\ldots, n-1 .
\eqno(2.22)
$$
In effect we are generating the polynomials orthogonal with respect
to the discrete inner product $(~ \cdot ~ , \cdot ~)_N$ in order to
approximate the desired orthogonal polynomials.

The computation of the approximate coefficients (2.22) can be done
by either\linebreak Stieltjes's procedure or Lanczos's algorithm (cf., e.g.,
[3, \S\S6--7]).
Both are implemented in the routine {\tt mcdis} of [4].

With any reasonable choice of the discretization (2.21), it will be
true that the procedure converges as $N \rightarrow \infty$,
$$
\lim_{N \rightarrow \infty} \hat{\alpha}_{k,N} = \hat{\alpha}_k , ~~
\lim_{N \rightarrow \infty} \hat{\beta}_{k,N} = \hat{\beta}_k , ~~~~
0 \leq k \leq n-1 .
\eqno(2.23)
$$
A natural choice, indeed, is given by
$$
\tau_k^{(N)} = t_k^{(N)} (d \lambda ), ~~~~
\omega_k^{(N)} = ~
\frac{w_k^{(N)} (d \lambda )}{\omega_m ( \tau_k^{(N)} )} ~, ~~~~
k = 1, 2,\ldots, N ,
\eqno(2.24)
$$
where $t_k^{(N)} (d \lambda )$ are the zeros of the orthogonal
polynomial $\pi_N (~\cdot~ ; d \lambda )$, and
$w_k^{(N)} (d \lambda )$ the respective Christoffel numbers.

The discretization method is conceptually simpler, and sometimes
more stable, than the methods of \S2.1, but may become significantly
more expensive, regardless of the choice of $m$, if poles are
close to the interval of integration,
or if high accuracy is desired; cf. Examples 3.1 and  3.5.
Note also that Case 4 that was skipped in \S2.1 can easily be
handled by the present method; see Example 3.6.
\bigskip

\noindent
{\bf 3. Numerical Examples}
\bigskip

All examples in this section were computed on the Cyber 205 in both
single and double precision.
The respective machine precisions are $7.11 \times 10^{-15}$ and
$5.05 \times 10^{-29}$.
\medskip

{\it Example} 3.1.  $I_1 (\omega ) = \int_{-1}^1 ~
\frac{ \pi t / \omega}{ \sin ( \pi t / \omega )} ~ dt , ~~~~
\omega > 1$.
\medskip

Here, $d \lambda (t) = dt$, and the poles of the integrand are
located at the integer multiples of $\omega$.
It is natural, then, to make our quadrature rule (1.8) exact for 
$m$ elementary rational functions matching
the $m$ poles closest to the
origin, say those at $- (m/2) \omega ,\ldots, - \omega$,
$\omega ,\ldots, (m/2) \omega$ when $m$ is even.
This suggests to identify $-1/ \zeta_\mu $ in (1.2) with these poles,
i.e., in (1.11a) to set 
$$
\xi_\nu = (-1)^\nu / ( \omega \lfloor (\nu +1)/2 \rfloor ),
~ \nu = 1, 2,\ldots, m. 
\eqno(3.1)
$$
Best accuracy is expected when $m = 2n$, in which case
the method described in \S2.1.1 was found to work rather well,
the only difficulty being the relatively slow convergence of
the backward recurrence algorithm for computing the 2$n$ modified
(Legendre) moments associated with the measure $dt/(t \pm \omega )$
when $\omega$ is very close to 1.
For single-precision accuracy
$\epsilon = \frac{1}{2}  \times 10^{-10}$ and double-precision
accuracy
$\epsilon^d  = \frac{1}{2} \times 10^{-25}$, and $n = 20$, the
respective starting indices $k_0$ and $k_0^d$ in the backward
recursion yielding the desired accuracy are shown in Table 3.1 for
selected values of $\omega$.

\bigskip
\begin{center}
\begin{tabular}{lrr}
\multicolumn{1}{c}{$\omega$} &
\multicolumn{1}{c}{$k_0$} &
\multicolumn{1}{c}{$k_0^d$} \\ \hline
\vrule height 2.5ex width0pt depth0pt\relax
2.0 & 50 & 63 \\
1.5 & 53 &71 \\
1.1 & 67 & 106 \\
1.01 & 124 & 247
\end{tabular}
\end{center}

\begin{center}
TABLE 3.1.  {\it Starting indices for backward \\
recurrence when} $n = 20$
\end{center}
\bigskip

Other than that, the method appears to be very stable and produces
quadrature rules that are rapidly converging.
In Table 3.2, the results of the $n$-point rule (1.17) 
\bigskip

\begin{center}
\begin{tabular}{lrlll}
\multicolumn{1}{c}{$\omega$} &
\multicolumn{1}{c}{$n$} &
\multicolumn{1}{c}{$n$-point rational Gauss} &
\multicolumn{1}{c}{$\gamma_n$} & 
\multicolumn{1}{c}{err. Gauss} \\ \hline
\vrule height 2.5ex width0pt depth0pt\relax
2.0 & 1 & 2.1 & 3.94(--1) & 1.43(--1) \\
& 4 & 2.33248722 & 3.50(--7) & 7.18(--5) \\
& 7 & 2.332487232246550235 & 2.61(--15) &2.73(--8) \\
& 10 & 2.332487232246550241107076 & 1.48(--24) & 1.02(--11) \\
1.1 & 2 & 4.43 & 1.73(--2) & 2.60(--1) \\
& 5 & 4.467773637 & 2.00(--9) & 2.09(--2) \\
& 8 & 4.46777364638776571 & 5.61(--18) & 1.53(--3)\\
& 11 & 4.467773646387765789236123 & 1.66(--27) & 1.09(--4)\\
1.01 & 3 & 8.429 & 2.53(--4) & 4.20(--1) \\
& 6 & 8.4301845803 & 6.27(--12) & 1.85(--1) \\
& 9 & 8.4301845804708420582 & 7.52(--21) & 8.37(--2) \\ 
& 12 & 8.430184580470842058971264 & 1.23(--30) & 3.75(--2)
\end{tabular}
\end{center}

\begin{center} TABLE 3.2.  {\it Numerical results for $I_1 ( \omega )$, error constants}, \\ 
{\it and comparison with Gauss quadrature}
\end{center}
\bigskip

\noindent
applied to $g(t) = (  \pi t / \omega ) / \sin ( \pi t/ \omega )$ in double
precision are shown for
$\omega = 2$, 1.1 and 1.01, along with the error constants
$\gamma_n$ of (1.18).
Also shown in the last column are the relative errors of the
$n$-point Gauss-Legendre rule.
For $\omega = 2$, the exact answer is known to be 8C/$\pi$, where C
is Catalan's constant (cf. [8, Eq. 3.747(2)]).
The value shown in Table 3.2 for $n = 10$ agrees with it to all 25
decimal digits given.
Ordinary Gauss-Legendre quadrature is seen
to converge rather slowly, as $\omega$ approaches 1.
In contrast, convergence of the rational Gauss quadrature rule is
fast even for $\omega$ very close to 1.
The extra effort required in this case is expended, as illustrated
in Table 3.1, at the time when the rule is generated.

\begin{center}
\begin{tabular}{lrlllr}
& & 
\multicolumn{2}{c}{method of \S2.1} &
\multicolumn{2}{c}{method of \S2.2} \\
\multicolumn{1}{c}{$\omega$} &
\multicolumn{1}{c}{$n$}\hspace{.25in} &
\multicolumn{1}{c}{SP} &
\multicolumn{1}{l}{\ \ DP} &
\multicolumn{1}{c}{SP} &
\multicolumn{1}{c}{\ \ \ \ \ DP} \\ \hline
\vrule height 2.5ex width0pt depth0pt\relax
2.0 & 1 \hspace{.25in} & .001 & .004 & .007 & .178 \\
& 4 \hspace{.25in} & .008 & .036 & .010 & .224 \\
& 7 \hspace{.25in} & .027 & .125 & .028 &.220 \\
& 10 \hspace{.25in} & .065 & .296 & .016 & .423 \\
1.1 & 2 \hspace{.25in} & .002 & .014 & .060 & 1.554  \\
& 5 \hspace{.25in} & .013 & .063 & .068 & 1.649 \\
& 8 \hspace{.25in} & .036 & .177 & .098 & 1.274 \\
& 11 \hspace{.25in} & .080 & .376 & .104 & 1.461 \\
1.01 & 3 \hspace{.25in} & .005 & .037 & .460 & 20.10 \,    \\
& 6 \hspace{.25in} & .020 & .104 & .446 & 41.51 \, \\
& 9 \hspace{.25in} & .050 & .244 & .458 & 10.44 \,  \\
& 12 \hspace{.25in} & .102 & .487 & .525 & 12.65 \,
\end{tabular}
\end{center}

\begin{center}
TABLE  3.3. {\it Timings} (in seconds) {\it of the methods in} \S\S2.1 \\
{\it and} 2.2 {\it applied to Example} 3.1
\end{center}

Identical results were obtained by the discretization method of
\S2.2, but with substantially greater effort, particularly for
higher accuracies and for $\omega$ close to 1.
Respective timings are shown in Table 3.3, both for
single-precision (SP) and double-precision (DP) accuracy
requirements of $\frac{1}{2} \times 10^{-10}$ and
$\frac{1}{2} \times 10^{-25}$, respectively.

While the choice $m = 2n$ indeed gives best accuracy, other choices of
$m$ may be preferable if the effort and time to generate the quadrature
rule is of any importance. 
It turns out that with the method of
partial fractions, $m = 2 \lfloor (n+1)/2 \rfloor$ gives almost the
same accuracy at about half the effort, whereas $m = 2$ gives
considerably less accuracy but requires only about one-tenth the
effort. The discretization method of \S2.2, on the other hand,
requires essentially the same effort regardless of the choice of $m$.
Some timings required to generate the quadrature rules for various
$m$ and $n$, and the relative errors achieved, are shown in Table 3.4
\bigskip

\begin{center}
\begin{tabular}{lcrccccl}
& & & \multicolumn{2}{c}{method of \S2.1} &
\multicolumn{2}{c}{method of \S2.2} \\
\multicolumn{1}{c}{$\omega$} &
\multicolumn{1}{c}{$n$} &
\multicolumn{1}{c}{$m$}\hspace{.25in} &
\multicolumn{1}{c}{SP} &
\multicolumn{1}{c}{DP} &
\multicolumn{1}{c}{SP} &
\multicolumn{1}{c}{DP} &
\multicolumn{1}{c}{err} \\ \hline
\\
2.0 & 10 & 2 \hspace{.25in} & .008 & .037 & .014 & \ \ .399 & 1.10(--17) \\
& & 10 \hspace{.25in} & .033 & .153 & .015 & \ \ .409 & 1.47(--25) \\
& & 20 \hspace{.25in} & .065 & .296 & .016 & \ \ .423 & 1.58(--25) \\
1.1 & 11 & 2 \hspace{.25in} & .009 & .047 & .095 & 1.402 & 2.20(--13) \\
& & 12 \hspace{.25in} & .045 & .213 & .100 & 1.426 & 2.80(--23) \\
& & 22 \hspace{.25in} & .080 &.376 & .104 & 1.461 & 6.68(--26) \\
1.01 & 12 & 2 \hspace{.25in} & .011 & .065 & .526 & 12.65 \ \ & 1.15(--13) \\
& & 12 \hspace{.25in} & .052 & .258 & .511 & 12.46 \ \ & 9.10(--25) \\
& & 24 \hspace{.25in} & .102 & .487 & .526 & 12.65 \ \ & 3.55(--27)
\end{tabular}
\end{center}

\begin{center}
TABLE 3.4.  {\it Timings and errors for selected} $m \leq 2n$
\end{center}
\bigskip

If $m = 2$, the values of $n$ for which full accuracy of about
$10^{-25}$ is attained are 15, 21 and 22 for $\omega = 2.0$, 1.1
and 1.01, respectively.
Interestingly, the timings involved are only about half those for
$m = 2 \lfloor (n+1) / 2 \rfloor$ shown in Table 3.4.

\bigskip
{\it Example} 3.2.
$I_2 ( \omega ) = \int_0^1 ~\frac{t^{- 1/2} \Gamma (1+t)}{t+
\omega} ~dt , ~~~~ 0 < \omega < 1$.
\medskip

Here we take $d \lambda (t) = t^{- 1/2}~ dt$ on [0,1]. If we wish
to match the first $2n-1$ poles
of the gamma function at the negative integers as well as the
pole at $- \omega$, we set $m = 2n$ in

\renewcommand{\arraystretch}{2}
$$
\begin{array}{l}
\xi_1 = \frac{1}{\omega} ~, \\
\xi_\nu = \frac{1}{\nu -1} , ~~~~ \nu = 2, 3 ,\ldots, m .
\end{array}
\eqno(3.2)
$$
\renewcommand{\arraystretch}{1}The rational
$n$-point Gauss rule (1.7), (1.17), generated by the method of \S2.1.1,
then produces (in double precision) results as shown in Table 3.5,
where $\omega = \frac{1}{2}$.
In the last column we list the absolute value of the difference
between double-precision and single-precision results.
In contrast to Example 3.1, we now see a case in which the accuracy
reaches a limit (at about $n = 10$) and deteriorates, rather than
improves, as $n$ is further increased.
(When $n = 20$, the calculation even breaks down in single
precision!).
The last column in Table 3.5 provides a clear hint as to what is
happening: a steady growth in numerical instability.
Closer examination reveals the true cause of this instability.
The constants $c_\nu$ in the partial fraction decomposition (2.7)
become very large and alternate in sign.
Thus, for example, $c_{18} = - 2.3375 \ldots \times 10^9$ and
$c_{19} = 2.3336 \ldots \times 10^9$ when $n = 18$.
This produces blocks of large coefficients $W_k$ in (2.10) that
alternate in sign from block to block, causing severe cancellation
errors in summations such as (2.5) (which are abundant in
Stieltjes's algorithm).
The phenomenon evidently is a manifestation of the asymmetric
distribution of the poles of the gamma function.
\bigskip

\begin{center}
\begin{tabular}{rll}
\multicolumn{1}{c}{$n$} &
\multicolumn{1}{c}{$n$-point rational Gauss} &
\multicolumn{1}{c}{$|DP - SP|$} \\ \hline
\vrule height 2.5ex width0pt depth0pt\relax
2 & 1.746 & 2.24(--13) \\
6 & 1.75012059121 & 7.05(--11) \\
10 & 1.750120591261335415386 & 3.36(--8) \\
14 & 1.7501205912613354159 & 1.52(--5) \\
18 & 1.7501205912613356 & 4.31(--3)
\end{tabular}
\end{center}

\begin{center}
TABLE 3.5.  {\it Numerical results for} $I_2 ( \omega )$, $\omega$ = .5
\end{center}
\bigskip

The method of \S2.2, in contrast, does not suffer from any
numerical instability and produces for $n = 12$, with comparable
effort, the value 
$$
I_2 (.5) = 1.750120591261335415394610,
\eqno(3.3)
$$
believed to be correct to all 25 digits shown.

Matching only $n$ poles, and thus taking $m = n$ in (3.2), stabilizes
the procedure of \S2.1.1 considerably, and as a consequence produces
the correct result (3.3) (except for a discrepancy of 1 unit in the
last decimal place) for $n = 11$. An even more stable procedure
results from taking $m = 2$ and, amazingly, yields the correct answer
(to all digits shown!) already for $n = 13$.
\bigskip

{\it Example} 3.3.
$I_3 ( \omega ) = \int_{-1}^1 \left( \frac{ \pi t/ \omega}{\sin (
\pi t/ \omega )} \right)^2 dt$, $~~\omega > 1$.
\medskip

Similarly as in Example 3.1, we take
$$
\xi_\nu = (-1)^\nu / ( \omega \lfloor ( \nu + 1 )/2 \rfloor ) ,
~~~ \nu = 1,2,\ldots, m/2 .
\eqno(3.4)
$$
We applied the procedure described in \S2.1.3 for $\omega = 2$, 1.5, 1.1
and 1.01, both in single and double precision, requesting accuracies
of $\epsilon = \frac{1}{2} \times 10^{-10}$ and $\epsilon^d = 
\frac{1}{2} \times 10^{-25}$, 
\bigskip

\begin{center}
\begin{tabular}{lllrr}
\multicolumn{1}{c}{$\omega$} &
\multicolumn{1}{c}{$k_1$} &
\multicolumn{1}{c}{$k_1^d$} &
\multicolumn{1}{c}{$k_2$} &
\multicolumn{1}{c}{$k_2^d$} \\ \hline
\vrule height 2.5ex width0pt depth0pt\relax
2.0 & 130 & 163 & 50 & 63 \\
1.5 & 133 & 191 & 53 & 71 \\
1.1 & 177 & 306 & 67 & 106 \\
1.01 & 384 & 727 & 124 & 247
\end{tabular}
\end{center}

\begin{center}
TABLE 3.6.  {\it Starting indices for two backward \\
recurrences when} $n = 20$
\end{center}

\noindent
respectively. When $m = 2n$ and $n = 20$, starting indices
$k_1$, $k_1^d$ in the first application
of the backward recursion that were found to meet the
accuracy requirements for the
poles closest to [--1,1], and the analogous starting indices $k_2$,
$k_2^d$ in the second application,
are shown in Table 3.6 for the four values of $\omega$. As expected,
the procedure becomes
laborious as $\omega$ approaches 1.
Selected double-precision results produced by the $n$-point rational Gauss
rule, along with error constants, are shown in Table 3.7.
The last column shows the relative error of results generated by the $n$-point
Gauss-Legendre rule.
For $\omega = 2$, the exact answer is known to be $I_3 (2) = 4$ ln
2 ([8, Eq. 3.837(2)]) and is correctly  reproduced to 25 digits
when $n = 11$.
Note again the fast convergence of the rational Gauss quadrature
rule, even for $\omega$ very close to 1, in contrast to the
relatively  slow convergence of the ordinary Gauss rule,
especially for $\omega$ close to 1.

\begin{center}
\begin{tabular}{lrlll}
\multicolumn{1}{c}{$\omega$} &
\multicolumn{1}{c}{$n$} &
\multicolumn{1}{c}{$n$-point rational Gauss} &
\multicolumn{1}{c}{$\gamma_n$} &
\multicolumn{1}{c}{err. Gauss} \\ \hline
\vrule height 2.5ex width0pt depth0pt\relax
2.0 & 2 & 2.75 & 9.90(--3) & 4.36(--2) \\
& 5 & 2.77258868 & 1.16(--9) & 4.70(--5) \\
& 8 & 2.7725887222397811 & 3.28(--18) & 2.92(--8) \\
& 11 & 2.772588722239781237668928 & 9.72(--28) & 1.52(--11) \\
1.1 & 2 & 15.5 & 3.00(--2) & 6.69(--1) \\
& 6 & 16.5328175 & 8.54(--12) & 6.79(--2) \\
& 10 & 16.5328177384604181 & 7.21(--24) & 3.42(--3) \\
& 14 & 16.53281773846041830155898 & 2.35(--37) & 1.40(--4) \\
1.01 & 2 & 184. & 6.34(--2) & 9.64(--1) \\
& 6 & 188.674782 & 2.20(--11) & 7.06(--1) \\
& 10 & 188.674784224994172 & 1.88(--23) & 4.00(--1) \\
& 14 & 188.6747842249941742708325 & 6.15(--37) & 1.92(--1)
\end{tabular}
\end{center}

\begin{center}
TABLE 3.7.  {\it Numerical results for} $I_3 (\omega )$, {\it error constants, and
comparison with Gauss quadrature}
\end{center}
\bigskip

While the choice $m = 2 \lfloor (n+1)/2 \rfloor$ produced similar
advantages as in Example 3.1 --- an increase of speed by a factor of
about 2 at only a slight loss of accuracy --- the choice $m = 2$
offered no significant gains in accuracy over the Gauss-Legendre rule,
unlike $m = 4$, which did (since a symmetric {\it pair} of double
poles is now accounted for).

We also applied the discretization method of \S2.2 and obtained
identical results with somewhat less effort in the case $\omega = 2$,
and about the same effort in the case $\omega = 1.1$.
For $\omega = 1.01$, however, we were unable to attain the
requested double-precision accuracy with a discretization parameter
$N \leq 800$ (in (2.21)).
\bigskip

{\it Example} 3.4.
$I_4 = \int_0^\infty ~\frac{t}{e^t -1}~e^{-t} dt$.
\medskip

The appropriate measure here is $d \lambda (t) = e^{-t} dt$ on
[0,$\infty$].
Since the integrand has poles at the integer multiples of
$2 \pi i$, we let $\zeta_\nu = -1 /(2 \nu \pi i ) = i/ (2 \nu \pi
)$, and thus in (1.12) take
$$
\xi_\nu = 0  , ~~~~
\eta_\nu = ~\frac{1}{2 \nu \pi} , ~~~~
\nu = 1,2,\ldots, m/2 .
\eqno(3.5)
$$
The quantity $p_\nu$ in (2.13) being real, and thus $d_\nu = 0$ in
(2.12), there is no nonconstant linear factor in the numerators of
(2.11). 
This simplifies somewhat the procedure in \S2.1.2, as it obviates
the need to apply the routine {\tt chri}.

In Table 3.8 we compare the performance (in double precision and for
$m = 2n$) of
our rational quadrature routine with Gauss-Laguerre quadrature
(applied to $f(t) = t/(e^t -1))$ and the Gaussian quadrature rule
(applied to $f(t) = e^{-t} )$ associated with ``Einstein's weight
function'' $t/(e^t - 1)$; for the latter see [6].
The respective relative errors are shown in the last two columns.
It can be seen that the Gauss-Laguerre and Gauss-Einstein
quadratures are
comparable in accuracy, the former being somewhat more accurate for
small values of $n$, the latter for larger values of $n$.
Both quadrature rules, however, are incomparably inferior to the
rational Gauss formula, which for $n = 15$ produces the true value
of the integral, $\zeta (2)-1  = (\pi^2 /6)-1$, to 25 correct
decimal digits.
(Actually, the last digit is off by one unit.) The results become even
slightly more accurate when we choose $m = 2 \lfloor (n+1)/2 \rfloor$,
and are still better, by several orders of magnitude, than those
for Gauss-Laguerre and Gauss-Einstein quadrature when $m = 2$.
\bigskip

\begin{center}
\begin{tabular}{rlll}
\multicolumn{1}{c}{$n$} &
\multicolumn{1}{c}{$n$-point rational Gauss} &
\multicolumn{1}{c}{err GL} &
\multicolumn{1}{c}{err GE} \\ \hline
\vrule height 2.5ex width0pt depth0pt\relax
1 & .59 & 9.76(--2) & 4.09(--1) \\
5 & .644934055 & 1.50(--5) & 2.97(--4) \\
10 & .644934066848226428 & 2.22(--8) & 1.15(--8) \\
15 & .6449340668482264364724151 & 1.59(--11) & 3.25(--13)
\end{tabular}
\end{center}

\begin{center}
TABLE 3.8.  {\it Numerical results for} $I_4$ {\it and comparison with \\
Gauss-Laguerre and Gauss-Einstein quadrature}
\end{center}

\bigskip

The high accuracy of our rational quadrature rules in this example
is all the more remarkable as the routine {\tt gchri}, used in their
construction (by the methods of \S2.1.2), is subject to ill-conditioning, causing the recursion
coefficients for the relevant orthogonal polynomials to gradually
lose accuracy (by as much as 10 decimals, when $n = 15$ and $m = 2n$).

This weakness is accentuated when one tries to deal with more
difficult integrals, for example,
$$
I (\theta ) = \int_0^\infty ~\frac{t}{e^t -1} ~
\sqrt{1 + \mbox{\small $\frac{1}{2}$} \theta t} ~dt = \int_0^\infty
,~~ \theta > 0,
\eqno(3.6)
$$
which has an additional
branch point singularity at $t = - 2 / \theta$.
Here, when $\theta = .75$, the $n$-point rational Gauss formula (in
double precision and for $m = 2n$) gives only about 13 correct
decimal places for $n = 15$, and 18 for $n = 30$.
By the time $n$ reaches 33, the ill-conditioning in the routine
{\tt gchri} has built up to such a level that the method fails (by
producing a negative $\beta$-recursion coefficient).
To get higher accuracy, one needs to apply the discretization
method of \S2.2, which is more stable, but becomes fairly expensive
if pushed much beyond $n  = 30$.
Using $d \lambda (t) = e^{-t} dt$, and hence the $N$-point
Gauss-Laguerre formula, to effect the discretization in (2.21), and
requesting an accuracy of $\frac{1}{2} \times 10^{-25}$ for the
desired recursion coefficients, we have observed timings of the
order 12--16 seconds, and discretization parameters $N$ as large as
$N = 370$, for $33 \leq n \leq 40$.
The rational Gauss formula so produced then yields relative errors
of $4.26 \times 10^{-21}$ for $n = 35$, and $9.65 \times 10^{-23}$
for $n = 40$.
This is still better, by about 4 decimal orders of accuracy, than
Gauss-Laguerre quadrature applied to the second form of the
integral in (3.6), and Gauss-Einstein quadrature applied to the
first form.

Generalized Fermi-Dirac integrals (cf. [12]) are similar to $I(
\theta )$ except that $t/(e^t -1)$ is replaced by $t^k /(e^{-\eta
+t} +1 )$, where $\eta$ is a real parameter and $k = 1/2$, 3/2 or
5/2.
The poles are now located at $\eta \pm (2 \nu -1) i \pi$, $\nu = 1,
2, 3,\ldots~$.
The use of rational Gauss quadrature to compute such integrals
is dealt with elsewhere [5].
\bigskip

{\it Example} 3.5.
$I_5 (\eta) = \int_0^{\infty} \frac{t}{e^{- \eta + t} - 1} e^{-t}
dt, ~~~ \eta < 0$.
\medskip

Again, we take $d \lambda (t) = e^{-t} dt$ and note that the poles are
now at $\eta \pm 2 \nu \pi i$, $\nu = 0, 1, 2, \ldots$ .
Accordingly, in (1.13) we take
$$
\xi_0  =  - \frac{1}{\eta}, ~~~
\xi_{\nu}  =  - \frac{\eta}{\eta^2 + 4 \nu^2 \pi^2}, ~~
\eta_{\nu} = \frac{2 \nu \pi}{\eta^2 + 4 \nu^2 \pi^2}, ~~~
\nu = 1, 2, \ldots, (m-1)/2 ,
\eqno{(3.7)}
$$
and use the procedure of \S2.1.$2^\prime$.  Selected results (for
$m = 2n-1$),
comparing rational Gauss formulae with Gauss-Laguerre formulae, are
shown in Table 3.9.  In the case 
\bigskip

\begin{center}
\begin{tabular}{rrll}
\multicolumn{1}{c}{$\eta$} &
\multicolumn{1}{c}{$n$} &
\multicolumn{1}{c}{$n$-point rational Gauss} &
\multicolumn{1}{c}{err GL} \\ \hline
\vrule height 2.5ex width0pt depth0pt\relax
-- .1  &  3 & .4503                          & 1.16(--1) \\
       &  6 & .4501936153                    & 5.13(--2) \\
       &  9 & .450193614441350               & 2.70(--2) \\
       & 12 & .45019361444134784096          & 1.55(--2) \\
--1.0  &  2 & .113                           & 2.14(--1) \\
       &  6 & .1111093520                    & 5.07(--3) \\
       & 11 & .1111093516052317322           & 1.81(--4) \\
       & 16 & .1111093516052317320105065     & 1.26(--5) \\
--10.0 &  2 & .122(--4)                       & 1.79(--1) \\
       &  6 & .113502121(--4)                & 1.57(--4) \\
       & 11 & .113502114635390578(--4)       & 7.31(--9) \\
       & 16 & .1135021146353905701870968(--4)& 1.20(--12)
\end{tabular}
\end{center}

\begin{center}
TABLE 3.9.  {\it Numerical results for} $I_5$ {\it and
comparison with Gauss-Laguerre quadrature}
\end{center}

\noindent
$\eta = - .1$, we were
able to go only up to $n = 13$; when $n = 14$, our procedure failed by
producing a negative $\beta$-coefficient in (2.2).  The difficulty
is caused by the ill-conditioning (mentioned after (2.10))
affecting the modified Chebyshev procedure.  Even though our
procedure was successful for $n = 13$, it had to work hard to take
care of the pole at $\eta = -.1$:  Backward recursion to compute
modified moments had to start at $\nu = 584$ to get
single-precision accuracy $\frac{1}{2} \times 10^{-10}$, and at
$\nu = 2650$ to get double-precision accuracy $\frac{1}{2} \times
10^{-25}$.

Replacing the numerator $te^{-t}$ in the integrand by $t^k$, where
$k = 1/2, 3/2$ or $5/2$, and adding a factor 
$\sqrt{1 + \frac{1}{2} \theta t}$ as in (3.6), produces the
Bose-Einstein integral whose computation by rational Gauss
quadrature is discussed in [5].
\bigskip

{\it Example} 3.6.
$I_6 = \int_0^\infty \left( \frac{t}{e^t -1} \right)^2 e^{-t} dt$.
\medskip

Here, as in Example 3.4, we take $d \lambda (t) = e^{-t} dt$ and
parameters $\xi_\nu$, $\eta_\nu$ as in (3.5), except that there are
only $m/4$ of them, $m$ being divisible by 4. In Table 3.10 we give
\bigskip

\begin{center}
\begin{tabular}{rlll}
\multicolumn{1}{c}{$n$} &
\multicolumn{1}{c}{$n$-point rational Gauss} &
\multicolumn{1}{c}{err GL} &
\multicolumn{1}{c}{err GE} \\ \hline
\vrule height 2.5ex width0pt depth0pt\relax
2 & .47 & 3.71(--2) & 1.61(--2) \\
8 & .4816405209 & 1.16(--6) & 5.99(--10) \\
14 & .4816405210580757311 & 4.36(--9) & 7.26(--18) \\
20 & .4816405210580757313458777 & 2.80(--11) & 1.09(--25)
\end{tabular}
\end{center}

\begin{center}
TABLE 3.10. {\it Numerical results for} $I_5$ {\it and comparison with
Gauss-Laguerre and Gauss-Einstein quadrature}
\end{center}
\bigskip

\noindent
the results for $m = 2n$ ($n$ even) obtained 
by the discretization method of \S2.2, analogous to those of Table 3.8
but using the square of the Einstein function as weight function in GE.
(The method of \S2.1, as mentioned earlier, was not implemented.)
What is remarkable in this example is the competitiveness of the
Gauss-Einstein quadrature rule vis-\`{a}-vis the rational Gauss rule.

\bigskip

\noindent
{\small {\bf References}}
\medskip

\begin{enumerate}
\itemsep=1pt
\parsep=1pt
\item 
{\tt GAUTSCHI, W.},
Minimal solutions of three-term recurrence relations and orthogonal
polynomials, {\it Math. Comp. 36} (1981), 547--554.
\item
{\tt GAUTSCHI, W.},
On generating orthogonal polynomials, {\it SIAM J. Sci. Stat.
Comput. 3} (1982), 289--317.
\item
{\tt GAUTSCHI, W.},
Computational problems and applications of orthogonal polynomials,
in {\it Orthogonal Polynomials and Their Applications} (C.
Brezinski et al., eds.),
IMACS Annals Comput. Appl. Math., Vol. 9, Baltzer, Basel, 1991, pp.
61--71.
\item 
{\tt GAUTSCHI, W.},
Algorithm xxx --- ORTHPOL:
A package of routines for generating orthogonal polynomials and
Gauss-type quadrature rules, {\it ACM Trans. Math. Software}, submitted.
\item
{\tt GAUTSCHI, W.}, On the computation of generalized Fermi-Dirac
and Bose-Einstein integrals, {\it Comput. Phys. Comm.}, to appear.
\item
{\tt GAUTSCHI, W.} and {\tt MILOVANOVI\'{C}, G.V.},
Gaussian quadrature involving Einstein and Fermi functions with an
application to summation of series, {\it Math. Comp. 44} (1985),
177--190.
\item 
{\tt GOLUB, G.H.} and {\tt WELSCH, J.H.},
Calculation of Gauss quadrature rules,
{\it Math. Comp. 23} (1969), 221--230.
\item 
{\tt GRADSHTEYN, I.S.} and {\tt  RYZHIK, I.M.},
{\it Table of Integrals, Series, and Products}, Academic Press,
Orlando, 1980.
\item 
{\tt L\'{O}PEZ LAGOMASINO, G.}, and {\tt ILL\'{A}N, J.},
A note on generalized quadrature formulas of Gauss-Jacobi type, in
{\it Constructive Theory of Functions}, Publ. House Bulgarian Acad.
Sci., Sofia, 1984, pp. 513--518.
\item 
{\tt L\'{O}PEZ LAGOMASINO, G.} and {\tt ILL\'{A}N GONZ\'{A}LEZ,
J.}, Sobre los m\'{e}todos interpolatorios de integraci\'{o}n
num\'{e}rica y su conexi\'{o}n con la aproximaci\'{o}n racional,
{\it Rev. Ciencias Mat\'{e}m. 8} (1987), no. 2, 31--44.
\item
{\tt PICHON, B.}, Numerical calculation of the generalized
Fermi-Dirac integrals, {\it Comput. Phys. Comm. 55} (1989),
127--136.
\item  
{\tt SAGAR, R.P.}, A Gaussian quadrature for the calculation of
generalized Fermi-Dirac integrals, {\it Comput. Phys. Comm. 66}
(1991), 271--275.
\item
{\tt VAN ASSCHE, W.}  and {\tt VANHERWEGEN, I.},
Quadrature formulas based on rational interpolation, {\it Math. Comp.},
to appear.
\end{enumerate}
\bigskip

\par
\indent
Professor Walter Gautschi

Department of Computer Sciences

Purdue University

West Lafayette, Indiana 47907, U.S.A.

wxg@cs.purdue.edu
\end{document}